\documentclass[12pt]{amsart}
\usepackage{amscd,amssymb,amsmath}
\usepackage[arrow,matrix]{xy}
\usepackage{graphicx}
\usepackage{cite}
\usepackage{ifpdf}

\newtheorem{thm}[subsection]{Theorem}
\newtheorem{lemma}[subsection]{Lemma}
\newtheorem{proposition}[subsection]{Proposition}

\numberwithin{equation}{section} 

\newcommand{\bea}{\begin{eqnarray}}
\newcommand{\eea}{\end{eqnarray}}

\begin{document}
\title [ The derivations of some evolution algebras]
{ The derivations of some evolution algebras}

\author {L. M. Camacho, J. R. G\'omez, B. A. Omirov, R. M. Turdibaev}

\address{L.\ M.\ Camacho and J. R. G\'omez\\ Dpto. Mathem\'atica Aplicada I, Universidad
de Sevilla, Avda. Reina Mercedes, s/n. 41012, Sevilla, Spain}
 \email {lcamacho@us.es, jrgomez@us.es}
 \address{B.\ A.\ Omirov\\ Institute of mathematics and information
 technologies, Do'rmon Yo'li str.,29, 100125,
Tashkent, Uzbekistan.} \email {omirovb@mail.ru}
 \address{R.\ M.\ Turdibaev\\ Department of Mathematics, National University of Uzbekistan, Vuzgorogok, 27, 100174,
Tashkent, Uzbekistan.} \email {rustamtm@yahoo.com}
 \maketitle
\begin{abstract}
In this work we investigate the derivations of  $n-$dimensional
complex evolution algebras, depending on the rank of the appropriate
matrices. For evolution algebra with non-singular matrices we prove
that the space of derivations is zero. The spaces of derivations for
evolution algebras with matrices of rank $n-1$ are described.

{\it AMS classifications (2010):} 17A36; 17D92; 17D99\\[2mm]

\end{abstract}
\maketitle
\section{Introduction and Preliminaries}

The notion of evolution algebras recently was introduced in the book
\cite{TIAN2}, where the author represented a lot of connections of
evolution algebras with the other objects in mathematics, genetic
and physics. The basic properties and some classes of evolution
algebras were studied as well in \cite{CLOR}, \cite{TIAN1},
\cite{TIAN2}.

 The concept of evolution algebras lies between
algebras and dynamical systems.  Although, evolution algebras do not
form a variety (they are not defined by identities), algebraically,
their structure has table of multiplication, which satisfies the
conditions of commutative Banach algebra. Dynamically, they
represent discrete dynamical systems. In fact, evolution algebras
are close related with graph and group theories, stochastic processes,
mathematical physics, genetics etc. The papers
\cite{ETH2}-\cite{GLI} were devoted to study of genetics using an abstract algebraic approach.

Recall the definition of evolution algebras.
 Let $E$ be a vector space over a field $K$ with defined multiplication
 $\cdot$ and a basis $\{e_1,e_2,\dots\}$ such that $$e_i\cdot e_j=0, \ i\neq j,$$
$$e_i\cdot e_i=\sum_k a_{ik}e_k, \ i \geq 1,$$ then $E$ is called evolution algebra and
basis $\{e_1,e_2,\dots\}$ is said to be natural basis.

From the above definition it follows that evolution algebras are commutative (therefore, flexible).

Let $E$ be a finite dimensional evolution algebra with natural
basis $\{e_1,\dots, e_n\},$ then
$$e_i\cdot e_i=\sum\limits_{j=1}^n a_{ij}e_j, \ 1\leq i \leq n,$$
where remaining products are equal to zero.

The matrix $A=(a_{ij})_{i,j=1}^{n}$ is called matrix of the algebra
$E$ in natural basis $\{e_1,\dots ,e_n\}.$

Obviously, $rank A=\dim (E\cdot E).$ Hence, for finite dimensional
evolution algebra the rank of the matrix does not depend on 
choice of natural basis.

The derivation for evolution algebra $E$ is defined as usual, i.e.,
a linear operator $d:E\to E$ is called a derivation if $$d(u\cdot
v)=d(u)\cdot v+u\cdot d(v)$$ for all $u,v \in E.$

Note that for any algebra, the space $Der(E)$ of all derivations is
a Lie algebra with the commutator multiplication.

Let $d$ be a derivation of evolution algebra $E$ with natural basis
$\{e_1,\dots,e_n\}$ and $d(e_i)=\sum_{j=1}^n d_{ij}e_j,\, 1\leq i
\leq n.$ Then the space of derivations for evolution algebra $E$ is
described as follows in \cite{TIAN2}.
$$Der(E)=\left\{ d\in End(E)\, |\, a_{kj}d_{ij}+a_{ki}d_{ij}=0, \textrm { for }
i\neq j; \, 2a_{ji}d_{ii}=\sum_{k=1}^n a_{ki}d_{jk} \right\}.$$

In the theory of non-associative algebras, particularly, in genetic algebras,
the Lie algebra of derivations of a given algebra is one of the important
tools for studying its structure. There has been much work on the subject of
derivations of genetic algebras (\cite{COSTA1},\cite{COSTA2},\cite{GONSHOR}).

In \cite{MICALI} multiplication is defined in terms of derivations,
showing the significance of derivation in genetic algebras. Several
genetic interpretations of derivation of genetic algebra are given
in \cite{HOLGATE}.

For evolution algebras the system of equations describing the derivations are
given in \cite{TIAN2}. In this work, we establish that the space of derivations
of evolution algebras with non-singular matrices is equal to zero.
The description of the derivations for evolution algebras,
the matrices of which are of rank $n-1$  is obtained.

\section{Main Result}
The following theorem describes derivations of evolution algebras with non-singular matrices.

\begin{thm} Let $d:E\to E$ be a derivation of evolution algebra $E$ with non-singular evolution matrix in basis $\langle e_1,\dots, e_n\rangle.$ Then this derivation $d$ is zero.
\end{thm}

\begin{proof} For a derivation  $d$ we have
$d(e_i)e_j+e_i d(e_j)=0$ and
$d(e_i e_i)=2d(e_i) e_i$ for all $1\leq i\neq j\leq n.$

Let $d(e_k)=\displaystyle  \sum_{i=1}^n d_{ki}e_i.$ Then we obtain $$d_{ij}(e_je_j)+d_{ji}(e_ie_i)=0 \eqno (1) $$
$$d(e_i e_i)=2d_{ii}(e_ie_i)\eqno (2) $$
for all $1\leq i\neq j\leq n.$

Since evolution matrix $A$ of algebra $E$ is non-degenerated and $(e_i\cdot e_i),\, (e_j\cdot e_j)$ represent the $i-$th and $j-$th rows of the matrix $A$ respectively, they can not be linearly dependent.

Thus, $d_{ij}=d_{ji}=0$ for all $1\leq i \neq j \leq n.$ Therefore, $d=diag\{d_{11},\dots, d_{nn}\}$ and $d(e_k)=d_{kk}e_k.$ Hence $spec(d)=\{d_{11},\dots, d_{nn}\}.$

Now $d(e_i\cdot e_i)=2d(e_i)\cdot e_i= 2d_{ii} (e_i\cdot e_i).$

Since $A$ is a non-singular, $e_i\cdot e_i \neq 0$ for all $1\leq i \leq n.$ The last equality shows that $$\{2d_{11},\dots, 2d_{nn}\} =spec(d).$$
This is possible if only $d$ is zero.
\end{proof}

Now we will investigate derivations for evolution algebras with matrices of rank $n-1.$

Since $rank A=n-1,$ making the suitable basis permutation we can assume that first $n-1$ rows of the matrix $A$ are linearly independent, i.e., $e_1e_1,\dots, e_{n-1}e_{n-1}$ are linearly independent and $e_ne_n=\displaystyle \sum_{k=1}^{n-1} b_k (e_ke_k)$ for some $b_1,\dots, b_{n-1} \in \mathbb{C}.$

Since $e_ie_i\neq 0$ for all $1\leq i\leq n-1,$  from $(2)$ we obtain that $2d_{ii}$ is an eigenvalue of $d$ for all $1\leq i\leq n-1.$ Hence,
$$spec(d)\supseteq \{2d_{11},2d_{22},\dots, 2d_{n-1n-1}\}.$$

Now from equality (1) we deduce $d_{ij}=d_{ji}=0$ for all $1\leq i\neq j\leq n-1.$

By putting $i=n$ to $(1)$ we obtain
$d_{nj}(e_je_j)+d_{jn}(e_ne_n)=0$ or
$$(d_{nj}+d_{jn}b_j)(e_je_j)+\sum_{k=1, k\neq j}^{n-1}d_{jn}b_k(e_ke_k)=0$$

Hence we obtain $d_{jn}b_k=0$ and $d_{nj}+d_{jn}b_j=0$ for all $1\leq k \neq j \leq n-1.$

Depending on different values of $b_k$ we will consider several cases.

\begin{lemma} Let $\displaystyle e_ne_n=\sum_{k=1}^{n-1} b_k (e_ke_k)$ and $b_p\neq 0,\,b_q\neq 0$ for some $1\leq p\neq q \leq n.$
Then $d=0.$ \end{lemma}

\begin{proof} In this case we have $d_{jn}b_p=0$ for all $1\leq j\neq p \leq n-1$ and $d_{jn}b_q=0$ for all $1\leq j\neq q \leq n-1,$ which implies $d_{jn}=0$ for all $1\leq j\leq n-1.$

Putting $d_{jn}=0$ to $d_{nj}+d_{jn}b_j=0$ we obtain $d_{nj}=0$ for all $1\leq j\leq n-1.$

Hence, $d=diag\{d_{11},d_{22},\dots, d_{nn}\}.$

Since $e_n e_n\neq 0$ from $(2)$ we obtain that $2d_{nn}$ is an eigenvalue of $d.$ Hence,
$$spec(d)=\{d_{11},d_{22},\dots, d_{nn}\}=\{2d_{11},2d_{22},\dots, 2d_{nn}\}$$ which is possible if only $d=0.$
\end{proof}

It should be noted that the opposite statement is not true.

From this lemma it follows that the only cases left to investigate
are $e_ne_n=0$ and $e_ne_n=b_ke_k,\, b_k\neq 0$ for some $1\leq k
\leq n.$ In the last case, by making suitable basis permutation one
can assume that $e_ne_n=b(e_1e_1), \, b\neq 0.$ Consider the
following $n\times n$ matrices:
  $$
\left(\begin{array}{cccccccc}
  d_{11} & 0 & \dots & 0 &0 & \dots &0 &  d_{1n}\\
   0     & 0 &\dots &0  & 0& \dots &0  &  0 \\
\vdots &    \vdots  & \ddots &  \vdots &  \vdots&& \vdots &  \vdots\\
 0     & 0 &\dots &0  & 0 &\dots &0  &  0 \\
       0     & 0 &\dots &0  & 2d_{11} &\dots &0  & 0\\
\vdots & \vdots & &\vdots & \vdots &\ddots & \vdots &\vdots\\
   0   & 0 &\dots & 0&  0 & \dots &2^{n-s-1}d_{11}&0  \\
   -bd_{1n}     & 0 &\dots &0  &0& \dots &0 &d_{11} \\
 \end{array}\right)\eqno (D_1) $$
$$
\left(\begin{array}{ccccccccccccc}
  d_{11} & 0 & \dots & 0 &0 & \dots&0  &0  &\dots&0&d_{1n}  \\
   0     & d_{22} &\dots &0  & 0&\dots& 0  &  0&\dots&0&0 \\
\vdots &\vdots &\ddots &\vdots &\vdots & &\vdots & \vdots&&\vdots&\vdots\\
     0     & 0 &\dots &2^{k-1}d_{22}  & 0 &\dots &0 &  0&\dots &0&0 \\
       0     & 0 &\dots &0  & 2d_{11} &\dots &0 & 0&\dots &0&0\\
      0   &0  & \dots &\vdots &\vdots & \ddots &  \vdots &\vdots&&\vdots&\vdots\\
   0   & 0 &\dots & 0&  0 & \dots &2^{m-k}d_{11}&0  &\dots &0&0\\
 0     & 0 &\dots &0  &0& \dots&0 &0 &\dots &0&0\\
   \vdots &\vdots &\ddots &\vdots &\vdots & &\vdots & \vdots&\ddots&\vdots&\vdots\\
   0     & 0 &\dots &0  &0& \dots&0 &0 &\dots &0&0\\
   -bd_{1n}     & 0 &\dots &0  &0& \dots&0 &0 &\dots &0&d_{11}\\
 \end{array}\right)\eqno (D_2)$$

 $$\left(\begin{array}{ccccccccc}
  d_{11} &0&\dots& 0 &  0 &\dots & 0&0&d_{1n}  \\
   0 &0&\dots& 0  & 0& \dots & 0&0&0  \\
   \vdots & \vdots &\ddots &\vdots  & \vdots &&\vdots&\vdots&\vdots\\
   0 &0&\dots& 0 &  0 &\dots & 0&0&0  \\
   0   &0&\dots & 0&\frac{d_{11}}{2^{n-s-2}} &\dots &0 &0 & 0 \\
   \vdots & \vdots &&\vdots  & \vdots &\ddots &\vdots&\vdots&\vdots\\
    0   & 0&\dots&0 & 0 &\dots&\frac{d_{11}}2&0  & 0 \\
   0   & 0&\dots&0 & 0 &\dots&0&d_{11}  & 0 \\
   -bd_{1n}  &0&\dots   & 0  &0&\dots  &0&0& d_{11} \\
 \end{array}\right) \eqno (D_3)$$

\begin{lemma}\label{lemma2} Let $e_n e_n =b (e_1 e_1), b\neq 0.$ Then derivation $d$ is either zero or it is in one of the following forms up to basis permutation:

$(i)$ $(D_1)$ where $d_{11}=\displaystyle \frac{\delta}{2^{n-s}-1},\, 1\leq s \leq n-1$ and $\delta^2=-bd_{1n}^2;$

$(ii)$ $(D_2)$ where  $d_{22}=\displaystyle\frac{1-2^{m-k}}{2^{k-1}}d_{11}, d_{11}=\displaystyle\frac{\delta}{2^{m-k+1}-1}, \, 1\leq k <m\leq n-1$ and $\delta^2=-bd_{1n}^2;$

$(iii)$ $(D_3)$ where $d_{11}=\delta$ and $\delta^2=-bd_{1n}^2.$
\end{lemma}

\begin{proof}
We have $d_{2n}=\dots=d_{n-1 n}=0,\, d_{n2}=\dots=d_{n n-1}=0$ and $d_{n1}=-bd_{1n}.$

By putting $i=n$ in $(2),$ we obtain $$2bd_{11}(e_1e_1)=b d(e_1e_1)=d(b(e_1e_1))=d(e_n e_n)=2d_{nn}(e_n e_n)=2d_{nn}b(e_1e_1).$$
Hence, $d_{11}=d_{nn}.$

From $(2)$ we deduce that
$$a_{i1} (d_{11}e_1+d_{1n}e_n)+ \sum_{j=2}^{n-1}a_{ij}d_{jj}e_j+
a_{in} (-bd_{1n}e_1+d_{11}e_n)= d(e_i e_i)=2d_{ii}(e_ie_i)=$$
$$2d_{ii}\sum_{j=1}^n a_{ij}e_j$$ which implies
$$a_{i1}(2d_{ii}-d_{11})=-a_{in}d_{1n}b \eqno (3)$$
$$a_{in}(2d_{ii}-d_{11})=a_{i1}d_{1n} \eqno (4)$$
$$a_{ij}(2d_{ii}-d_{jj})=0 \eqno (5)$$ for all $1\leq i \leq n-1$ and $2\leq j\leq n-1.$

If $d_{1n}=0,$ then $d=diag\{d_{11},\dots, d_{n-1n-1},d_{11}\}$ and $\{d_{11},\dots, d_{n-1n-1}\}=spec(d)\supseteq\{2d_{11},2d_{22},\dots, 2d_{n-1n-1}\}$ which leads to $d=0.$

Assume that $d_{1n}\neq 0.$ One can find $spec(d)=\{d_{22},\dots,d_{n-1 n-1}, \alpha ,\beta \},$  where $\alpha=d_{11}+\delta ,\, \beta=d_{11}- \delta$ and $\delta^2=-b d_{1n}^2.$ Obviously, $\alpha \neq \beta.$

Let $\lambda \in spec(d)$ be such that $|\lambda|=\max \{|\alpha|,|\beta|,|d_{22}|,\dots, |d_{n-1n-1}|\}.$

If $\lambda \in \{d_{22},\dots, d_{n-1 n-1}\}$ then $2\lambda$ is also an eigenvalue which contradicts to module maximality of $\lambda.$ Therefore $\lambda =\alpha$ or $\lambda =\beta.$

Also note that from $(3)$ and $(4)$ it follows that $a_{i1}=0$ if and only if $a_{in}=0.$

If $a_{i1} \neq 0 \,(a_{in} \neq 0),$ then multiplying $(3)$ and $(4)$ we obtain $(2d_{ii}-d_{11})^2=-bd_{1n}^2$ or $2d_{ii}=d_{11}\pm \delta.$ Hence for these $i$ we have $$d_{ii}=\frac12 \alpha \textrm{ or } d_{ii}=\frac12\beta.\eqno(6)$$

Now we consider several cases.\\

\textbf{\textit{Case 1.}} Let $\alpha \beta \neq 0, \alpha+\beta \neq 0.$
Since $\alpha+\beta = 2d_{11}\in spec(d)$ and $\alpha+\beta \not\in\{\alpha,\beta\}$ we obtain that there exists $i_1$ such that $d_{i_1i_1}= \alpha+\beta.$ Then $2d_{i_1 i_1} \in spec(d)$ which implies that $2d_{i_1i_1}=d_{i_2 i_2}$ for some $i_2$ or $2d_{i_1i_1}\in \{\alpha, \beta\}.$ If $2d_{i_1i_1}=d_{i_2 i_2}$ we can continue till we obtain $2^k d_{i_1 i_1}=\dots = 2d_{i_k i_k}\in \{\alpha, \beta\}$ for some $1\leq k\leq n-2.$

Thus, for some $1\leq k\leq n-2$ we have $2^k(\alpha+\beta)\in \{ \alpha, \beta\}.$

Let us assume that $2^k (\alpha+\beta)=\alpha.$

Then $d_{11}=\frac{\alpha}{2^{k+1}}, \, d_{i_1i_1}=\frac{\alpha}{2^{k}},\dots, d_{i_k i_k}=\frac{\alpha}{2}$ and   $\beta=-(1-\frac1{2^k})\alpha.$ Hence, $|\beta|<|\alpha|$ and obviously, $2^s \beta \neq 2^r \alpha$ for any $r,s\in Z.$

Consider the possible non-zero values of $|d_{22}|,\dots,|d_{n-1
n-1}|$ and let them be $d_1<\dots<d_p.$ We already know that
$\{d_1,\dots, d_{n-1}\}\supseteq
\{\frac{|\alpha|}{2^k},\dots,\frac{|\alpha|}{2}\}.$ Since
$spec(d)\supseteq \{2d_{22},\dots,2d_{n-1 n-1}\}$ one obtains that
$2d_1,\dots, 2d_p \in \{d_1,\dots,d_p,|\alpha|,|\beta|\}.$

Since $2d_p \leq |\alpha|$ and $d_{i_ki_k}=\frac{\alpha}2$ we conclude that  $d_p=\frac{|\alpha|}{2}.$

Observe that there  can be only one eigenvalue $d_{i_ki_k}=\frac{\alpha}2$ with module $d_p.$ Indeed, if for some $i$ we have $|d_{ii}|=d_p, d_{ii}\neq \frac{\alpha}2,$ then $spec(d)\ni 2d_{ii}\neq \alpha$ and $|2d_{ii}|=|\alpha|.$ Therefore, there exists $j$ such that $d_{jj}=2d_{ii}.$ But then $2d_{jj}\in spec(d)$ and $|2d_{jj}|=2|\alpha|>|\alpha|$ which is a contradiction.

Now since there is only one eigenvalue with module $\frac12 |\alpha|$ one obtains that there is only one eigenvalue $\frac14 \alpha$ of module $d_{p-1}$ and etc.

If not all $d_2,\dots,d_p$ are in the form $\frac1{2^m}|\alpha|$ then applying similar arguments to $|\beta|$ we obtain that there can be at most only one eigenvalue $\frac12 \beta$ with module $\frac12|\beta|,$ $\frac14\beta$ with module $\frac14 |\beta|$ and etc.

Hence, $\displaystyle \{d_{22},\dots ,d_{n-1 n-1}\} \setminus \{0\}=\bigcup_{i=1}^s \{ \frac1{2^i} \alpha\}$ or $\displaystyle \{d_{22},\dots ,d_{n-1 n-1}\} \setminus \{0\}=\bigcup_{i=1}^s \{ \frac1{2^i} \alpha\} \bigcup_{j=1}^r \{ \frac1{2^j} \beta\}.$

\textbf{\textit{Case 1.1.}} Let $\displaystyle \{d_{22},\dots ,d_{n-1 n-1}\} \setminus \{0\}=\bigcup_{i=1}^s \{ \frac1{2^i} \alpha\}.$ Then from (6) for those $i$ such that $a_{in}\neq 0$ we obtain $2d_{ii}=\alpha.$ Then $(4)$ implies $a_{i1}=\frac{\alpha-d_{11}}{d_{1n}}a_{in}.$ Hence, the first and the last columns of the matrix $A$ are collinear. Therefore, all other columns must be non-zero and linearly independent so that $rank A=n-1$ is satisfied.

Assume that there are $s-1$ zeros among $d_{22},\dots,d_{n-1n-1}.$ Then $0=d_{22}=\dots =d_{ss}<|d_{s+1 s+1}|\leq |d_{s+2 s+2}|\leq \dots \leq |d_{n-1n-1}|.$
If $2\leq i\leq s$ and $s+1\leq j\leq n-1$ then $2d_{ii}-d_{jj}=-d_{jj}\neq 0$ and from $a_{ij}(2d_{ii}-d_{jj})=0$ we obtain $a_{ij}=0$ for $2\leq i\leq s,\,s+1\leq j\leq n-1.$

Now if $2\leq j\leq s$ and $s+1\leq i\leq n-1$ then $d_{jj}=0,d_{ii}\neq 0$ and from $a_{ij}(2d_{ii}-d_{jj})=0$ we conclude that $a_{ij}=0$ for $2\leq j\leq s,\,s+1\leq i\leq n-1.$

Since $d_{s+1s+1}\neq 2d_{ii}$ for all $2\leq i \leq n-1,$ from $(5)$ we obtain $a_{is+1}=0$ for all $2\leq i \leq n-1.$ Since $2-,\dots,(n-1)-$th columns are linearly independent, $a_{1s+1}\neq 0$ and therefore, by $(5)$ we obtain $d_{s+1s+1}=2d_{11}.$

Now we will show that among $d_{s+1 s+1},\dots,d_{n-1n-1}$ there are no equal elements. Let $d_{s+1 s+1}=d_{s+2 s+2}.$ Then $d_{s+2s+2}\neq 2d_{ii}$ for all $2\leq i \leq n-1$ and from $(5)$ we obtain $a_{is+2}=0$ for all $2\leq i \leq n-1.$ Hence, the $(s+2)-$th column is either zero or collinear to $(s+1)-$th column of matrix $A.$ This is a contradiction.

Now let $|d_{s+1 s+1}|<|d_{s+2 s+2}|<\dots<|d_{jj}|=|d_{j+1 j+1}|\leq \dots \leq |d_{n-1n-1}|$ for some $s+2\leq j\leq n-2.$ Then $2d_{ii}-d_{jj}=0$ if and only if $i=j-1$ and therefore $(5)$ implies $a_{ij}=0$ for all $i\neq j-1.$ Similarly, since $2d_{ii}-d_{j+1j+1}=0$ only if $i=j-1$ we obtain $a_{ij+1}=0$ for all $i\neq j-1.$ This implies that either columns $j$ and $j+1$ are collinear or at least one of them is zero, which is a contradiction.
Hence in this case all $d_{s+1 s+1},\dots ,d_{n-1n-1}$ are distinct and $2^{n-s} d_{11}=2^{n-s-1}d_{s+1s+1}=\dots=2d_{n-1n-1}=\alpha$ and hence $2^{n-s}d_{11}=d_{11}+ \delta \Rightarrow d_{11}=\frac1{2^{n-s}-1}\delta.$

Therefore matrix $A$ should be in the form
$$\left(
\begin{array}{ccccllccc}
0 &0      &\dots & 0  & a_{1s+1} & 0& \dots   & 0& 0  \\
0 &a_{22} &\dots & a_{2s} & 0 & 0 & \dots  & 0& 0  \\
\vdots &\vdots &      & \vdots & \vdots &\vdots  &\ddots &  \vdots   &\vdots  \\
0 &a_{s2} &\dots &a_{ss}  & 0 & 0 & \dots   & 0 & 0  \\
0 & 0     &\dots & 0      & 0 & a_{s+1 s+2}  &\dots & 0 & 0  \\
\vdots &   \vdots   &   & \vdots& \vdots &\vdots&\ddots & \vdots &\vdots  \\
0 &0 &\dots & 0           & 0& 0 & \dots  & a_{n-2 n-1} & 0  \\
a_{n-1 1} &0 & \dots & 0           & 0 & 0& \dots  & 0 & a_{n-1n}\\
0  & 0&\dots & 0           & ba_{1s+1} & 0& \dots  & 0 & 0\\
\end{array}\right) \eqno (A_1)$$
and $d$ is in the form $(D_1)$ with  $d_{11}=\displaystyle \frac{\delta}{2^{n-s}-1}.$

\textbf{\textit{Case 1.2.}} Let $\displaystyle \{d_{22},\dots ,d_{n-1 n-1}\} \setminus \{0\}=\bigcup_{i=1}^s \{ \frac1{2^i} \alpha\} \bigcup_{j=1}^r \{ \frac1{2^j} \beta\}.$

Assume that  $\{d_{22},\dots,d_{kk}\}=\displaystyle \bigcup_{j=1}^r \{ \frac1{2^j} \beta\}, \, \{d_{k+1 k+1},\dots,d_{mm}\}= \displaystyle\bigcup_{i=1}^s \{ \frac1{2^i} \alpha\}$ and $d_{m+1m+1}=\dots =d_{n-1n-1}=0$ such that $|d_{22}| \leq \dots \leq |d_{kk}|,\, |d_{k+1k+1}| \leq \dots \leq |d_{mm}|.$

Since $2d_{ii}-d_{22}\neq 0$ for all $1\leq i \leq n$ from $(5)$ and due to $a_{in}=ba_{i1}$ we obtain $a_{i2}=0$ for all $1\leq i \leq n.$ Now since $rank A =n-1,$ the other columns must be non-zero and linearly independent. Similarly, as in Case 1.1 one obtains that $d_{33}\neq d_{22}$ and so on.

Hence, $$d_{kk}=2d_{k-1 k-1}=\dots=2^{k-1}d_{22}$$ and for all $3\leq j \leq k$ it follows that $a_{j-1j}\neq 0, a_{ij}=0\,(i\neq j-1).$

Now since $2d_{ii}-d_{k+1k+1}\neq 0$ for all $2\leq i \leq n-1,$ where $d_{k+1k+1}$ is $\frac1{2^s}\alpha,$ it must be $d_{k+1k+1}=2d_{11}.$ Otherwise, the $(k+1)-$ column is zero, which is a contradiction. Then in the $(k+1)-$th column the only non-zero elements are $a_{1k+1}$ and $a_{nk+1}=ba_{1k+1}.$

Applying the similar arguments as in Case 1.1 we deduce $$d_{mm}=2d_{m-1m-1}=\dots=2^{m-k-1}d_{k+1k+1}=2^{m-k}d_{11}$$ and for all $k+1\leq j \leq m$ we have $a_{j-1j}\neq 0, a_{ij}=0, i\neq j-1.$

Now for all $1\leq i\leq m$ and $m+1\leq j \leq n-1$ we have $2d_{ii}-d_{jj}=2d_{ii}\neq 0.$ Then from $(6)$ we obtain $a_{ij}=0$ for all $1\leq i\leq m$ and $m+1\leq j \leq n-1.$ Also, since $a_{nj}=ba_{1j},$ it follows $a_{nj}=0$ for $m+1\leq j \leq n.$

Hence, $d_{kk}=\frac12 \beta$ and $d_{mm}=\frac12\alpha.$ Then from $(4)$ and $(6)$ it follows that\\ $\displaystyle a_{kn}=\frac{d_{1n}}{\beta-d_{11}}a_{k1}\neq 0$ and $\displaystyle
a_{mn}=\frac{d_{1n}}{\alpha-d_{11}}a_{m1}\neq 0.$

Also from $d_{11}+\delta= \alpha=2d_{mm}=2^{m-k+1}d_{11}$  it follows that  $d_{11}=\displaystyle\frac{\delta}{2^{m-k+1}-1}.$

Now $2d_{11}-\alpha=\beta=2d_{kk}=2^kd_{22}$ implies $d_{22}=\displaystyle\frac{1-2^{m-k}}{2^{k-1}}d_{11}.$

Hence, the matrix of $A$ is
$$\left(
\begin{array}{ccccccccccccc}
0&0 &0 &\dots &0 &a_{1k+1} &0 & \dots& 0&0 &\dots &0&0\\
0&0 &a_{23} &\dots &0 &0 & 0& \dots& 0 & 0 &\dots &0&0\\
\vdots&\vdots &\vdots &\ddots &\vdots &\vdots & \vdots& & \vdots & \vdots & &\vdots&\vdots\\
0& 0&0 &\dots &a_{k-1k} & 0&0 &\dots &0 & 0 &\dots &0&0\\
a_{k1}&0 &0 & \dots& 0 &0 &0& \dots &0 & 0 &\dots &0&a_{kn}\\
0& 0&0 & \dots& 0 & 0& a_{k+1k+2} & \dots& 0& 0 &\dots &0&0\\
\vdots&\vdots &\vdots & &\vdots &\vdots &\vdots &\ddots & \vdots& \vdots & &\vdots& \vdots\\
0& 0&0 & \dots& 0 & 0& 0 & \dots& a_{m-1m}& 0 &\dots &0&0\\
a_{m1}&0 &0 & \dots &0 &0 &0 & \dots&0 & 0 &\dots &0&a_{mn}\\
0& 0& 0& \dots &0 &0 & 0& \dots&0& a_{m+1 m+1}&\dots &a_{m+1 n-1}&0\\
\vdots&\vdots &\vdots & &\vdots & \vdots&\vdots & \dots&\vdots & \vdots & \ddots&\vdots&\vdots\\
0& 0& 0& \dots &0 &0 & 0& \dots&0&a_{n-1 m+1}&\dots &a_{n-1 n-1}&0\\
0&0 &0 & \dots &0 & ba_{1k}& 0& \dots&0 &0 &\dots &0&0\\
\end{array}\right) \eqno $$

Denote by $(A_2)$ the form of the above matrix. For the evolution algebra with matrix in the form $(A_2)$ the derivation $d$ is in the form $(D_2)$ with   $d_{22}=\displaystyle\frac{1-2^{m-k}}{2^{k-1}}d_{11}$ and
 $d_{11}=\displaystyle\frac{\delta}{2^{m-k+1}-1} .$\\

Note that, we can assume $2^k(\alpha+\beta)=\beta$ in the beginning of our argumentation in this case. Then in Case 1.1 we obtain that $d$ is in the form $(D_1)$ with  $d_{11}=\displaystyle \frac{-\delta}{2^{n-s}-1}.$ In Case 1.2 $d$ is in the form $(D_2)$ with $d_{22}=\displaystyle\frac{1-2^{m-k}}{2^{k-1}}d_{11}$ and
 $d_{11}=\displaystyle\frac{-\delta}{2^{m-k+1}-1}.$\\

\textbf{\textit{Case 2.}} Let $\alpha \beta \neq 0, \alpha=-\beta,$ i.e., $d_{11}=0.$ We will show that this case is impossible.

Obviously, there are non-zero elements among $d_{22},\dots, d_{n-1n-1}.$ Otherwise, from $(3)$ and $(4)$ it follows that the first and the last columns of matrix $A$ are zero, which is a contradiction.

Now consider the possible non-zero values of $|d_{22}|,\dots,|d_{n-1 n-1}|$ and let them be $d_1<\dots<d_p.$  Since $spec(d)\supseteq \{2d_{22},\dots,2d_{n-1 n-1}\}$ one obtains that $2d_1,\dots, 2d_p \in \{d_1,\dots,d_p,|\alpha|\}.$

Since this values are non-zero, we deduce that  $|\alpha|=2d_p,d_p=2d_{p-1},\dots,d_2=2d_1.$

Observe that there can be only eigenvalue $\frac{\alpha}2$ or $-\frac{\alpha}2$ with module $d_p.$ Indeed, if for some $i$ we have $|d_{ii}|=d_p, d_{ii}\neq \pm\frac{\alpha}2$ we obtain $spec(d)\ni 2d_{ii}\neq \pm\alpha$ and $|2d_{ii}|=|\alpha|.$ Therefore, there exists $j$ such that $d_{jj}=2d_{ii}.$ But then $2d_{jj}\in spec(d)$ and $|2d_{jj}|=2|\alpha|>|\alpha|$ which is a contradiction.

Now since $\pm \frac12 \alpha$ are the only possible eigenvalues with module $\frac12 |\alpha|$ one obtains that the only possible eigenvalues with module  $d_{p-1}$ are $\pm\frac14 \alpha$ and etc.

Hence,  $\displaystyle \{d_{22},\dots ,d_{n-1 n-1}\} \setminus \{0\}\subseteq\bigcup_{i=1}^s \{ \frac1{2^i} \alpha\} \bigcup_{j=1}^r \{ -\frac1{2^j} \alpha\}.$

If $\frac12 \alpha \not\in \{d_{22},\dots,d_{n-1 n-1}\}$ and $-\frac12 \alpha \not\in \{d_{22},\dots,d_{n-1 n-1}\}$ then from $(6)$ we obtain that the first and the last columns are zero which contradicts to $rank A=n-1.$ Hence there exists $2\leq k\leq n-1$ such that $d_{kk}\in \{\frac12\alpha, -\frac12\alpha\}.$

Now, if $\displaystyle \{d_{22},\dots ,d_{n-1 n-1}\} \setminus \{0\}\supseteq \bigcup_{i=1}^s \{ \frac1{2^i} \alpha\}$ and $-\frac12 \alpha \not \in \{d_{22},\dots ,d_{n-1 n-1}\}$ then by $(4)$ and $(6)$ we obtain that the first and the last columns of matrix $A$ are linearly dependent, i.e.,  $a_{i1}=\frac{\alpha}{d_{1n}}a_{in}$ for all $1\leq i \leq n.$ Hence, in order to be $rank A=n-1$ the other columns must be non-zero and linearly independent.

However, if $d_{pp}=\frac1{2^s}\alpha,$ then from $(5)$ we obtain that the $p-$th column is zero which is a contradiction.

Now if $\displaystyle \{d_{22},\dots ,d_{n-1 n-1}\} \setminus \{0\}\supseteq \bigcup_{j=1}^r \{ -\frac1{2^j} \alpha\}$ and $\frac12 \alpha \not \in \{d_{22},\dots ,d_{n-1 n-1}\}$ then by $(4)$ and $(6)$ we obtain that the first and the last columns of matrix $A$ are linearly dependent, i.e.,  $a_{i1}=\frac{-\alpha}{d_{1n}}a_{in}$ for all $1\leq i \leq n.$ Hence, in order to be $rank A=n-1$ the other columns must be non-zero and linearly independent.

However, if $d_{pp}=\frac1{2^q}\alpha,$ then from $(5)$ we obtain that the $p-$th column is zero which is a contradiction.

Now let $\displaystyle \{d_{22},\dots ,d_{n-1 n-1}\} \setminus \{0\}=\bigcup_{i=1}^s \{ \frac1{2^i} \alpha\} \bigcup_{j=1}^r \{ -\frac1{2^j} \alpha\}.$
Then for $d_{pp}=\frac1{2^s} \alpha$ and $d_{qq}=-\frac1{2^r} \alpha$ we obtain $2d_{ii}-d_{pp}\neq 0, 2d_{ii}-d_{qq}\neq 0$ for all $1\leq i \leq n-1$ and hence from $(5)$ the $p-$th and $q-$th columns are zero which is a contradiction to $rank A=n-1.$\\

\textbf{ Case 3.} Let $\alpha \neq 0,\beta =0.$

Then $2d_{11}=\alpha=d_{11}+\delta,$ and hence $d_{11}=\delta.$

Let us consider the possible non-zero values of $|d_{22}|,\dots,|d_{n-1 n-1}|$ and let them be $d_1<\dots<d_p.$  Since $spec(d)\supseteq \{2d_{22},\dots,2d_{n-1 n-1}\}$ one obtains that $2d_1,\dots, 2d_p \in \{d_1,\dots,d_p,|\alpha|\}.$

Since this values are non-zero, it follows that $|\alpha|=2d_p,d_p=2d_{p-1},\dots,d_2=2d_1.$

Observe that there can be only eigenvalue $\frac{\alpha}2$ with module $d_p.$ Indeed,if for some $i$ we have $|d_{ii}|=d_p, d_{ii}\neq \frac{\alpha}2$ we obtain $spec(d)\ni 2d_{ii}\neq \alpha$ and $|2d_{ii}|=|\alpha|.$ Therefore, there exists $j$ such that $d_{jj}=2d_{ii}.$ But then $2d_{jj}\in spec(d)$ and $|2d_{jj}|=2|\alpha|>|\alpha|$ which is a contradiction.

Similarly, since there is only one eigenvalue with module $\frac12 \alpha$ one obtains that there is only one eigenvalue $\frac14 \alpha$ of module $d_{p-1}$ and etc.

Thus, $spec(d)=\{\frac1{2^p}\alpha,\dots,\frac12\alpha,\alpha\}$ or $spec(d)=\{0,\frac1{2^p}\alpha,\dots,\frac12\alpha,\alpha\}.$ Again, by making suitable basis permutation one can assume that $|d_{22}|\leq\dots\leq |d_{n-1n-1}|.$

Assume that there are $s-1$ zeros among $d_{22},\dots,d_{n-1n-1}.$ Then $0=d_{22}=\dots =d_{ss}<d_{s+1 s+1}\leq d_{s+2 s+2}\leq \dots \leq d_{n-1n-1}.$
If $1\leq i\leq s$ and $s+1\leq j\leq n-1$ then $2d_{ii}-d_{jj}\neq 0$ and from $a_{ij}(2d_{ii}-d_{jj})=0$ we obtain $a_{ij}=0$ for $1\leq i\leq s,\,s+1\leq j\leq n-1.$

Now if $2\leq j\leq s$ and $s+1\leq i\leq n-1$ then $d_{jj}=0,d_{ii}\neq 0$ and from $a_{ij}(2d_{ii}-d_{jj})=0$ we obtain $a_{ij}=0$ for $2\leq j\leq s,\,s+1\leq i\leq n-1.$

Since $d_{s+1s+1}\neq 2d_{ii}$ for all $s+1\leq i \leq n-1,$ from $a_{is+1}(2d_{ii}-d_{s+1s+1})=0$ we obtain $a_{is+1}=0$ for all $s+1\leq i \leq n-1,$ i.e., the $(s+1)-$th column of matrix $A$ is zero.

Now we will show that among $d_{s+1 s+1},\dots,d_{n-1n-1}$ there are no equal elements. Let $d_{s+1 s+1}=d_{s+2 s+2}.$ Then $d_{s+2s+2}\neq 2d_{ii}$ for all $s+1\leq i \leq n-1,$ from $a_{is+2}(2d_{ii}-d_{s+2s+2})=0$ we obtain $a_{is+2}=0$ for all $s+1\leq i \leq n-1$ i.e., the $(s+2)-$th column of matrix $A$ is zero. This is a contradiction to $rank A=n-1.$

Now let $d_{s+1 s+1}<d_{s+2 s+2}<\dots<d_{jj}=d_{j+1 j+1}\leq \dots \leq d_{n-1n-1}$ for some $s+2\leq j\leq n-2.$ Then $2d_{ii}-d_{jj}=0$ only if $i=j-1$ and therefore $a_{ij}(2d_{ii}-d_{jj})=0$ implies $a_{ij}=0$ for all $i\neq j-1.$ Similarly, since $2d_{ii}-d_{j+1j+1}=0$ only if $i=j-1$ we obtain $a_{ij+1}=0$ for all $i\neq j-1.$ This implies that either columns $j$ and $j+1$ are collinear or at least one of them is zero. However, this contradicts to $rank A=n-1.$
Hence in this case all $d_{s+1 s+1},\dots ,d_{nn}$ are distinct.

Also from $(6)$ it follows that $a_{i1}=a_{in}=0$ for all $s+1\leq i \leq n-1.$

Therefore matrix $A$ should be in the form
$$\left(
\begin{array}{ccccccccc}
a_{11}&0 & \dots & 0 & 0  & 0& \dots  & 0& a_{1n}  \\
a_{21}&a_{22} & \dots & a_{2s} & 0& 0 & \dots   & 0& a_{2n}  \\
\vdots&\vdots &       & \vdots & \vdots &\vdots  & & \vdots    &\vdots  \\
a_{s1}&a_{s2} & \dots & a_{ss} & 0 & 0 & \dots  & 0 & a_{sn}  \\
0&0 & \dots & 0           & 0 & a_{s+1 s+2} &\dots& 0  & 0  \\
\vdots&\vdots &       & \vdots & \vdots &\vdots&\ddots   &  \vdots   &\vdots  \\
0&0 & \dots & 0  & 0  & 0& \dots & a_{n-2 n-1} & 0  \\
a_{n-1 1}&0 & \dots & 0           & 0  & 0& \dots & 0 & a_{n-1n}\\
ba_{1n}  &0& \dots & 0 & 0 &0& \dots & 0 & 0 \\ \end{array}\right). \eqno(A_3)$$

Hence, for the evolution algebra with matrix in the form $(A_3)$ the derivation $d$ is in the form $(D_3)$ with $d_{11}=\delta.$

Note that in symmetrical case $\alpha=0, \beta \neq 0$  one can obtain in a similar way that $d$ is in the form $(D_3)$ with $d_{11}=-\delta.$  So the statement of Lemma \ref{lemma2} is verified.
\end{proof}

The following lemma completes the description of derivations of evolution algebras with matrices of rank $n-1.$

\begin{lemma} Let evolution algebra has a matrix $A=(a_{ij})_{1\leq i,j\leq n}$
in the natural basis $e_1,\dots,e_n$ such that $e_n e_n=0$ and $rank A=n-1.$
Then derivation $d$ of this evolution algebra is
either zero or it is in one of the following forms up to basis permutation:

$$\left(\begin{array}{cccc}
0&\dots&0&d_{1n}\\
\vdots&\ddots&\vdots&\vdots\\
0&\dots&0&d_{n-1n}\\
0&\dots&0&0\\
\end{array}\right),\eqno (D_4)$$
 where $\displaystyle \sum_{k=1}^{n-1} a_{ik}d_{kn}=0,$ $1\leq i \leq n-1;$

$$\left(
\begin{array}{ccccccc}
0&\dots&0&0&\dots&0&0\\
\vdots&\ddots&\vdots&\vdots&&\vdots&\vdots\\
0&\dots&0&0&\dots&0&0\\
0&\dots&0&\frac{d_{nn}}{2^{n-k-1}}&\dots&0&d_{k+1n}\\
\vdots&\ddots&\vdots&\vdots&\ddots&\vdots&\vdots\\
0&\dots&0&0&\dots&\frac{d_{nn}}2&d_{n-1n}\\
0&\dots&0&0&\dots&0&d_{nn}\\
\end{array}\right),  \eqno (D_5)$$
 where $d_{i+1 n}=\frac{a_{in}}{a_{ii+1}}\left(\frac1{2^{n-i-1}}-1\right)d_{nn},\,a_{ii+1}\neq0,$ $k+1\leq i \leq n-2, 1\leq k \leq n-1$ and $d_{k+1n}\in\mathbb{C}.$
\end{lemma}

\begin{proof} From $e_n e_n=0$ we obtain  $d_{nj}=0$ for all $1\leq j\leq n-1.$ Now one can see that
$spec(d)=\{d_{11},\dots,d_{nn}\}\supseteq\{2d_{11},2d_{22},\dots, 2d_{n-1n-1}\}.$

Let $\lambda \in spec(d)$ be such that $|\lambda|=\max_{1\leq i \leq n} |d_{ii}|.$

If $\lambda \in \{d_{11},\dots,d_{n-1n-1}\} $ then $2\lambda \in spec(d)$ which yields $\lambda =0.$ Therefore, in this case we obtain $d_{11}=\dots=d_{nn}=0$ and $d(e_i)=d_{in}e_n$ for all $1\leq i \leq n-1.$ Then from $(2)$ it follows that $\displaystyle \sum_{j=1}^n a_{ij}d_{jn} e_n=d(e_i e_i)=0$ for all $1\leq i \leq n-1.$ The last one implies that vector $(d_{1n},\dots,d_{n-1 n},0)$ is a solution of homogeneous linear system of equations $Ax=0.$ Observe that if the first $n-1$ columns are linearly independent then $d=0.$

In order to $d\neq 0$ we consider the matrices with first $n-1$ columns linearly dependent. Denote the form of this matrices by $(A_4).$

So in this case $d$ is in the form $(D_4).$\\

Now if $\lambda \notin \{d_{11},\dots,d_{n-1n-1}\}$ then $\lambda =d_{nn}$ and we can assume that $d_{nn} \neq 0.$ Consider the possible non-zero values of $|d_{11}|,\dots,|d_{n-1 n-1}|$ and let them be $d_1<\dots<d_p.$ Since $spec(d)\supseteq \{2d_{11},\dots,2d_{n-1 n-1}\}$ one obtains that $2d_1,\dots, 2d_p \in \{d_1,\dots,d_p,|d_{nn}|\}.$ Since this values are non-zero, we deduce $|d_{nn}|=2d_p,d_p=2d_{p-1},\dots,d_2=2d_1.$ Observe that there  can be only one eigenvalue $\frac12 d_{nn}$ with module $d_p.$ Indeed,if for some $i<n$ we have $|d_{ii}|=d_p, d_{ii}\neq d_{nn}$ we obtain $spec(d)\ni 2d_{ii}\neq d_{nn}$ and $|2d_{ii}|=|d_{nn}|.$ Therefore, there exists $1\leq j \leq n-1$ such that $d_{jj}=2d_{ii}.$ But then $2d_{jj}\in spec(d)$ and $|2d_{jj}|=2|d_{nn}|>|d_{nn}|$ which is a contradiction. Similarly, since there is only one eigenvalue with module $\frac12 d_{nn}$ one obtains that there is only one eigenvalue $\frac14 d_{nn}$ of module $d_{p-1}$ and etc.

Hence, $spec(d)=\{\frac{d_{nn}}{2^p},\dots,\frac{d_{nn}}2,{d_{nn}}\}$ or $spec(d)=\{0,\frac{d_{nn}}{2^p},\dots,\frac{d_{nn}}2,{d_{nn}}\}.$ Now making appropriate basis permutation we can assume that $|d_{11}|\leq\dots\leq |d_{n-1n-1}|<|d_{nn}|.$

From $(2)$ we obtain $$\sum_{j=1}^{n-1}a_{ij}d_{jj}e_j+ \sum_{j=1}^n(a_{ij}d_{jn})e_n=d(e_i e_i)=2d_{ii}(e_ie_i)=2d_{ii}\sum_{j=1}^n a_{ij}e_j,$$ which implies $\displaystyle \sum_{j=1}^na_{ij}d_{jn}=2d_{ii}a_{in}$ and $a_{ij}(2d_{ii}-d_{jj})=0$ for all $1\leq i,j\leq n-1.$

Assume that there are $k$ zeros among $d_{11},\dots,d_{n-1n-1}.$ Then $0=d_{11}=\dots =d_{kk}<|d_{k+1 k+1}|\leq \dots \leq |d_{n-1n-1}|<|d_{nn}|.$
If $1\leq i\leq k$ and $k+1\leq j\leq n-1$ then $d_{ii}=0,d_{jj}\neq 0$ and from $a_{ij}(2d_{ii}-d_{jj})=0$ it follows that $a_{ij}=0$ for $1\leq i\leq k,\,k+1\leq j\leq n-1.$

Analogously, if $1\leq j\leq k$ and $k+1\leq i\leq n-1$ then $d_{jj}=0,d_{ii}\neq 0$ and from $a_{ij}(2d_{ii}-d_{jj})=0$ we obtain $a_{ij}=0$ for $1\leq j\leq k,\,k+1\leq i\leq n-1.$

Since $d_{k+1k+1}\neq 2d_{ii}$ for all $k+1\leq i \leq n-1,$ from $a_{ik+1}(2d_{ii}-d_{k+1k+1})=0$ we obtain $a_{ik+1}=0$ for all $k+1\leq i \leq n-1,$ i.e., the $(k+1)-$th column of matrix $A$ is zero.

Now we will show that among $d_{k+1 k+1},\dots,d_{nn}$ there are no equal elements. Let $d_{k+1 k+1}=d_{k+2 k+2}.$ Then $d_{k+2k+2}\neq 2d_{ii}$ for all $k+1\leq i \leq n-1,$ from $a_{ik+2}(2d_{ii}-d_{k+2k+2})=0$ we obtain $a_{ik+2}=0$ for all $k+1\leq i \leq n-1$ i.e., the $(k+2)-$th column of matrix $A$ is zero. This is a contradiction to $rank A=n-1.$

Now let $|d_{k+1 k+1}|<|d_{k+2 k+2}|<\dots<|d_{jj}|=|d_{j+1 j+1}|\leq \dots <|d_{nn}|$ for some $k+2\leq j\leq n-2.$ Then $2d_{ii}-d_{jj}=0$ only if $i=j-1$ and therefore $a_{ij}(2d_{ii}-d_{jj})=0$ implies $a_{ij}=0$ for all $i\neq j-1.$ Similarly, since $2d_{ii}-d_{j+1j+1}=0$ only if $i=j-1$ we obtain $a_{ij+1}=0$ for all $i\neq j-1.$ This implies that either columns $j$ and $j+1$ are collinear or at least one of them is zero. However, this contradicts to $rank A=n-1.$
Hence in this case all $d_{k+1 k+1},\dots ,d_{nn}$ are distinct and $d_{ii}=\frac{d_{nn}}{2^{n-i}}$ for all $k+1\leq i \leq n-1.$

Now if $k+1\leq i,j\leq n-1$ we have $2d_{ii}-d_{jj}=0$ if and only if $j=i+1$ and hence $a_{ij}=0$ for all $k+1\leq i \leq n-1$ and $k+1\leq j \leq n-1, j\neq i+1.$ Therefore, matrix $A$ should be in the form
$$\left(
\begin{array}{ccccccccc}
a_{11} & \dots & a_{1k} & 0 &0& \dots &  0& a_{1n}  \\
\vdots & \ddots& \vdots & \vdots &\vdots  &      &\vdots &\vdots \\
a_{k1} & \dots & a_{kk} & 0 & 0&\dots & 0   & a_{kn}  \\
0 & \dots & 0           & 0 & a_{k+1 k+2} & \dots &0 &a_{k+1 n}  \\
\vdots &       & \vdots & \vdots &\vdots&\ddots       &\vdots &\vdots \\
0 & \dots & 0           & 0 & 0&\dots &  a_{n-2 n-1} & a_{n-2n}  \\
  0 & \dots & 0           & 0 &0& \dots  & 0 & a_{n-1n}\\
  0 & \dots & 0           & 0 &0& \dots  & 0 & 0\\
  \end{array}\right) \eqno(A_5)$$

Denote by $A_k=(a_{ij})_{1\leq i,j \leq k}$ the $k\times k$ submatrix  of matrix $A.$

Since $rank A =n-1$ we obtain $\det  A_k \cdot a_{k+1 k+2}\cdot\dots\cdot a_{n-1 n} \neq 0.$

Now $\displaystyle \sum_{j=1}^na_{ij}d_{jn}=2d_{ii}a_{in}$ implies
$$\left(
\begin{array}{ccc}
a_{11} & \dots & a_{1k}  \\
\vdots &       & \vdots  \\
a_{k1} & \dots & a_{kk} \\ \end{array}\right)
\left(
\begin{array}{c}
d_{1n}\\
\vdots\\
d_{kn}\\
\end{array}\right) =
\left(
\begin{array}{c}
0\\
\vdots\\
0\\
\end{array}\right)\eqno (7)$$
and  $a_{i i+1}d_{i+1n}+a_{in}d_{nn}=2d_{ii}a_{in}$ for all $k+1\leq i \leq n-2$ and $a_{n-1n}d_{nn}=2d_{n-1 n-1}a_{n-1n}$ which is an identity.

Now since $\det A_k\neq 0$ from $(7)$ it follows that $d_{1n}=\dots=d_{kn}=0.$

For $k+1\leq i \leq n-2$ we obtain $a_{in}(2d_{ii}-d_{nn})=a_{ii+1}d_{i+1n}$ which implies $$ d_{i+1 n}=\frac{a_{in}}{a_{ii+1}}\left(\frac1{2^{n-i-1}}-1\right)d_{nn}.$$ Hence, derivation $d$ is in the form of $(D_5).$
\end{proof}

As a result of previous lemmas we obtain the following

\begin{thm} Let $d:E\to E$ be a derivation of $n-$dimensional evolution algebra $E$ with matrix $A$ in basis $\langle e_1,\dots, e_n\rangle$ such that $rank A=n-1.$ Then the derivation $d$ is either zero or is in one of the forms given in Lemma 1.2 and Lemma 1.3.
\end{thm}

We can conclude that if the matrix of evolution algebra $E$ can be transformed by basis permutation to matrices of the form $(A_1)-(A_5),$ then in this permuted basis the corresponding derivations are in the form $(D_1)-(D_5),$ respectively.
Moreover, if the matrix of evolution algebra $E$ can not be transformed by basis permutation to any of the forms $A_i, 1\leq i \leq 5,$ then derivation of such algebra is zero.

For all $1\leq i \leq 5$ denote by $E_i$ an evolution algebra with matrix, that can be transformed by basis permutation to the form $A_i.$

Then it is easy to see that $\dim Der(E_i)=2, i\neq 4$ and $\dim Der(E_4)=1.$

\begin{proposition} Let evolution algebra $E_{(k)} (1\leq k \leq n)$ with natural basis $\{e_1,\dots, e_n\}$ be such that $\displaystyle e_ie_i =\sum_{j=i}^k a_{ij}e_j, a_{ii}\neq 0$ for $1\leq i \leq k$ and $e_i e_i=0$ for $k+1\leq i \leq n.$ Then in this basis the derivation has the following matrix: $$\left(\begin{array}{cc}
O & O\\
O & D\\
\end{array}\right)\eqno (8)$$ where $D\in M_{n-k}(\mathbb{C}).$
\end{proposition}

\begin{proof} From $(1)$ it follows that $d_{ij}(e_je_j)+d_{ji}(e_ie_i)=0$ for all $1\leq i\neq j \leq n.$ Now if we take $1\leq i\neq j\leq k$ then $e_ie_i$ and $e_je_j$ are linearly independent. Hence, we obtain $d_{ij}=d_{ji}=0$ for all $1\leq i\neq j \leq k.$

Now if $1\leq i \leq k$ and $k+1\leq j \leq n$ then $e_je_j=0$ and hence $d_{ji}(e_ie_i)=0.$ This implies that $d_{ji}=0$ for all $1\leq i \leq k,k+1\leq j \leq n.$

From $(2)$ we have $d(e_ke_k)=2d(e_k)e_k.$

Since $\displaystyle d(e_ke_k)=a_{kk}d(e_k)=a_{kk}\sum_{j=k}^n d_{kj}e_j$ and
$2d(e_k)e_k=2d_{kk}(e_ke_k)=2d_{kk}a_{kk}e_k$ we obtain
$d_{kk}=d_{kk+1}=\dots=d_{kn}=0.$

Assume that $d(e_{k-j+1})=\dots=d(e_k)=0$ for some $j.$

From $(2)$ we have $d(e_{k-j}e_{k-j})=2d(e_{k-j})e_{k-j}.$

Since $\displaystyle d(e_{k-j}e_{k-j})= d(\sum_{p=k-j}^k a_{{k-j}p}e_{p})=a_{k-j\,k-j}d(e_{k-j})=a_{k-j\,k-j}\sum_{p=k-j}^n d_{k-jp}e_p$ and
$\displaystyle 2d(e_{k-j})e_{k-j}=2d_{k-j k-j}(e_{k-j}e_{k-j})=2d_{k-j k-j}\sum_{p=k-j}^k a_{k-jp}e_p$ we deduce
$d_{k-jk-j}=0$ and hence $d_{k-j\,k-j+1}=\dots=d_{k-jn}=0.$
Therefore, $d(e_{k-j})=0$ and we obtain $d(e_1)=\dots=d(e_k)=0.$

Since for $k+1 \leq i,j\leq n$ equalities $(1)$ and $(2)$ turn into identities, we obtain that $d$ is in the form $(8).$
\end{proof}

Note that $\dim Der(E_{(k)})=(n-k)^2.$\\

\textbf{Acknowledgments.} This work is supported in part by the PAICYT, FQM143 of Junta de Andaluc\'ia (Spain). The third named author was supported by the grant NATO-Reintegration ref. CBP.EAP.RIG. 983169. The last named author would like to acknowledge ICTP OEA-AC-84 for a given support.

\end{document}